\documentclass[12pt]{article}
\usepackage{amsmath,amsfonts,amssymb,amsthm}
\input amssym.def
\begin{document}
\def \Z{\Bbb Z}
\def \C{\Bbb C}
\def \R{\Bbb R}
\def \Q{\Bbb Q}
\def \N{\Bbb N}

\def \A{{\mathcal{A}}}
\def \D{{\mathcal{D}}}
\def \E{{\mathcal{E}}}
\def \E{{\mathcal{E}}}
\def \H{\mathcal{H}}
\def \S{{\mathcal{S}}}
\def \wt{{\rm wt}}
\def \tr{{\rm tr}}
\def \span{{\rm span}}
\def \Res{{\rm Res}}
\def \Der{{\rm Der}}
\def \End{{\rm End}}
\def \Ind {{\rm Ind}}
\def \Irr {{\rm Irr}}
\def \Aut{{\rm Aut}}
\def \GL{{\rm GL}}
\def \Hom{{\rm Hom}}
\def \mod{{\rm mod}}
\def \ann{{\rm Ann}}
\def \ad{{\rm ad}}
\def \rank{{\rm rank}\;}
\def \<{\langle}
\def \>{\rangle}

\def \g{{\frak{g}}}
\def \h{{\hbar}}
\def \k{{\frak{k}}}
\def \sl{{\frak{sl}}}
\def \gl{{\frak{gl}}}

\def \be{\begin{equation}\label}
\def \ee{\end{equation}}
\def \bex{\begin{example}\label}
\def \eex{\end{example}}
\def \bl{\begin{lem}\label}
\def \el{\end{lem}}
\def \bt{\begin{thm}\label}
\def \et{\end{thm}}
\def \bp{\begin{prop}\label}
\def \ep{\end{prop}}
\def \br{\begin{rem}\label}
\def \er{\end{rem}}
\def \bc{\begin{coro}\label}
\def \ec{\end{coro}}
\def \bd{\begin{de}\label}
\def \ed{\end{de}}

\newcommand{\m}{\bf m}
\newcommand{\n}{\bf n}
\newcommand{\nno}{\nonumber}
\newcommand{\nord}{\mbox{\scriptsize ${\circ\atop\circ}$}}
\newtheorem{thm}{Theorem}[section]
\newtheorem{prop}[thm]{Proposition}
\newtheorem{coro}[thm]{Corollary}
\newtheorem{conj}[thm]{Conjecture}
\newtheorem{example}[thm]{Example}
\newtheorem{lem}[thm]{Lemma}
\newtheorem{rem}[thm]{Remark}
\newtheorem{de}[thm]{Definition}
\newtheorem{hy}[thm]{Hypothesis}
\makeatletter \@addtoreset{equation}{section}
\def\theequation{\thesection.\arabic{equation}}
\makeatother \makeatletter

\begin{center}
{\Large \bf Quantum vertex $\C((t))$-algebras and quantum affine
algebras}
\end{center}

\begin{center}
{Haisheng Li\footnote{Partially supported by NSF grant
DMS-0600189}\\
Department of Mathematical Sciences\\
Rutgers University, Camden, NJ
08102}
\end{center}

\begin{abstract}
We give a summary of the theory of (weak) quantum vertex
$\C((t))$-algebras and the association of quantum affine algebras
with (weak) quantum vertex $\C((t))$-algebras.
\end{abstract}

\section{Introduction}
In the earliest days of vertex (operator) algebra theory, Lie
algebras had played an important role, and in particular, an
important family of vertex operator algebras (see \cite{flm},
\cite{fz}, \cite{dl}) was associated with untwisted affine Lie
algebras. A fundamental problem, posed in \cite{fj} (cf.
\cite{efk}), has been to establish a suitable theory of quantum
vertex algebras, so that quantum affine algebras can be canonically
associated with quantum vertex algebras in the same (or a similar)
way that affine Lie algebras are associated with vertex operator
algebras. In literature, there have been a few of notions of quantum
vertex (operator) algebra (\cite{efr}, \cite{ek}, \cite{b-qva},
\cite{li-qva1}, \cite{ab}), however this particular problem is still
to be solved.

In a series of papers, starting with \cite{li-qva1}, we have been
investigating vertex algebra-like structures arising from various
algebras including quantum affine algebras and Yangians, with an
ultimate goal to solve the aforementioned problem. Our key idea is
to start with the algebraic structures that the generating functions
of those quantum algebras on highest weight modules could possibly
``generate.'' This is the fundamental guideline of this series of
studies.

The first paper \cite{li-qva1} was to provide a foundation for the
whole series. Starting with an {\em arbitrary} vector space $W$, we
studied general (formal) vertex operators (=quantum fields) on $W$,
which are elements of $\Hom (W,W((x)))$. Let $\E(W)$ denote the
space $\Hom (W,W((x)))$ alternatively. Then we studied various types
of subsets of $\E(W)$ and the algebraic structures generated by such
subsets, where the operations of vertex operators on vertex
operators are given by what is often called ``operator product
expansion.'' The most general type consists of what were called
``quasi compatible subsets,'' whereas compatible subsets are
relatively restrictive, but still very general. It was proved
therein that any (quasi) compatible subset of $\E(W)$ generates a
nonlocal vertex algebra with $W$ as a (quasi) module in a certain
sense. (Nonlocal vertex algebras are analogs of noncommutative
associative algebras, in contrast to that vertex algebras are
analogs of commutative and associative algebras.) This generalizes
the main result of \cite{li-g1}, which states that every compatible
subset generates a nonlocal vertex algebra with $W$ as a module. It
follows from this general result that a wide variety of algebras can
be associated with nonlocal vertex algebras. In particular, if $W$
is taken to be a highest weight module for a quantum affine algebra,
the generating functions in the Drinfeld realization form a quasi
compatible subset of $\E(W)$, and therefore they generate a nonlocal
vertex algebra with $W$ as a quasi module.

Furthermore, with the defining relations of quantum affine algebras
in mind, we formulated and studied a notion of ``pseudo local
subset'' of $\E(W)$ with $W$ a general vector space, to single out a
family of quasi compatible subsets. Now, given a pseudo local subset
$U$, one has a nonlocal vertex algebra $\<U\>$ with $W$ as a quasi
module. Under a certain assumption, we proved in \cite{li-qva1} that
there exists a unitary quantum Yang-Baxter operator on $\<U\>$ with
two spectral parameters, which describes the braided commutativity
relation of the vertex operators from the set $\<U\>$. Roughly
speaking, one obtains a deformed chiral algebra structure in the
sense of \cite{efr}. Note that this quantum Yang-Baxter operator is
for vertex operators on {\em the quasi module} $W$, not for the
adjoint vertex operators on {\em the algebra} $\<U\>$.

Motivated by Etingof-Kazhdan's notion of quantum vertex operator
algebra \cite{ek}, in particular by the $\S$-locality axiom, we
studied in \cite{li-qva1} a notion of ``$\S$-local subset'' of
$\E(W)$ (with $W$ a vector space), which singles out a family of
compatible subsets. It was proved that if $U$ is an $\S$-local
subset of $\E(W)$, the adjoint vertex operators on the nonlocal
vertex algebra $\<U\>$ satisfies $\S$-locality (commutativity). This
lead us to a theory of (weak) quantum vertex algebras and their
modules. A conceptual result is that for any $\S$-local subset $U$
of $\E(W)$, $\<U\>$ is a weak quantum vertex algebra with $W$ as a
canonical module. As the set of the Drinfeld generating functions of
quantum affine algebras is not $\S$-local, this theory of (weak)
quantum vertex algebras leaves quantum affine algebras out.
Nevertheless, it has been proved to be suitable for studying
Yangians. More specifically, in \cite{li-infinity} we have
associated certain versions of double Yangians with quantum vertex
algebras. Furthermore, in \cite{li-hqva}, we formulated a notion of
$\hbar$-adic (weak) quantum vertex algebra and we associated
centrally extended double Yangians with $\hbar$-adic quantum vertex
algebras and their modules.

In \cite{li-tqva}, we came back to the problem with quantum affine
algebras again. On the basis of \cite{li-qva1}, we developed a
theory of (weak) quantum vertex $\C((t))$-algebras and we
successfully associated quantum affine algebras with weak quantum
vertex $\C((t))$-algebras. In this theory, a weak quantum vertex
$\C((t))$-algebra is a $\C((t))$-module and a nonlocal vertex
algebra over $\C$, that satisfies a Jacobi-like identity. As a
nonlocal vertex algebra, a weak quantum vertex $\C((t))$-algebra
satisfies the associativity for ordinary vertex algebras.
Furthermore, a quantum vertex $\C((t))$-algebra is a weak quantum
vertex $\C((t))$-algebra equipped with a unitary quantum Yang-Baxter
operator on $V$ with two formal parameters, which describes the
braided commutativity relation of vertex operators $Y(v,x)$ for
$v\in V$ and satisfies some other properties. Even though a weak
quantum vertex $\C((t))$-algebra $V$ is a $\C((t))$-module, the
vertex operator map $Y$ is not $\C((t))$-linear, as by definition
$$Y(f(t)u,x)g(t)v=f(t+x)g(t)Y(u,x)v\ \ \ \mbox{
for }f(t),g(t)\in \C((t)),\; u,v\in V$$ (where linearity is
deformed). Thus,  the formal variable $t$ in this theory is not a
deformation parameter, unlike the formal variable $\hbar$ in
Etingof-Kazhdan's theory of quantum vertex operator algebras
\cite{ek}.

The notion of weak quantum vertex $\C((t))$-algebra naturally arisen
{}from our study on the nonlocal vertex algebras generated by pseudo
local subsets of $\E(W)$ with $W$ a vector space over $\C$. (Recall
that any quasi compatible subset $U$ of $\E(W)$ generates a nonlocal
vertex algebra $\<U\>$ over $\C$.)  It has been realized in
\cite{li-qva1} that $\<U\>$ is not large enough to describe the
braided commutativity relation and one needs to consider the span
$\C((x))\<U\>$, noticing that $\E(W)\ (=\Hom (W,W((x))))$ is
naturally a $\C((x))$-module. It has also been proved therein that
$\C((x))\<U\>$ is a nonlocal vertex algebra over $\C$, but it is not
a nonlocal vertex algebra over $\C((x))$, as the adjoint vertex
operator map is not $\C((x))$-linear. This lead us to a notion of
nonlocal vertex $\C((t))$-algebra, where a nonlocal vertex
$\C((t))$-algebra $V$ is simply a $\C((t))$-module and a nonlocal
vertex algebra over $\C$ whose vertex operator map $Y$ satisfies the
deformed $\C((t))$-linear property mentioned previously. In terms of
this notion, for any quasi compatible subset $U$ of $\E(W)$,
$\C((x))\<U\>$ is a nonlocal vertex $\C((t))$-algebra with $f(t)\in
\C((t))$ acting as $f(x)$. Furthermore, we considered what we called
``quasi $\S(x_{1},x_{2})$-local subsets'' and
``$\S(x_{1},x_{2})$-local subsets'' of $\E(W)$, where quasi
$\S(x_{1},x_{2})$-local subsets are quasi compatible while
$\S(x_{1},x_{2})$-local subsets are compatible. (The notion of quasi
$\S(x_{1},x_{2})$-local subset is a slight reformulation of the
notion of pseudo local subset.) Our conceptual result is that for
any (quasi) $\S(x_{1},x_{2})$-local subset $U$ of $\E(W)$,
$\C((x))\<U\>$ is a weak quantum vertex $\C((t))$-algebra with
$f(t)\in \C((t))$ acting as $f(x)$, and the space $W$ becomes a
(quasi) module. Furthermore, to construct a quantum vertex
$\C((t))$-algebra from a weak quantum vertex $\C((t))$-algebra, we
extended Etingof-Kazhdan's notion of non-degeneracy for nonlocal
vertex $\C((t))$-algebras and we proved that just as with quantum
vertex algebras \cite{li-qva2}, every non-degenerate weak quantum
$\C((t))$-algebra has a (unique) canonical quantum vertex
$\C((t))$-algebra structure.

Now, taking $W$ to be a highest weight module for a quantum affine
algebra and taking $U$ to be the set of the Drinfeld generating
functions, one has a quasi $\S(x_{1},x_{2})$-local subset $U$ of
$\E(W)$, and then one has a weak quantum vertex $\C((t))$-algebra
$\C((x))\<U\>$ with $W$ as a canonical quasi module. To a certain
extent, this solves the problem mentioned at the very beginning,
though we still have to show that this weak quantum vertex
$\C((t))$-algebra is a quantum vertex $\C((t))$-algebra, or
sufficiently to show that it is non-degenerate.

We mention that there is a very interesting and closely related work
\cite{ab}, in which Anguelova and Bergvelt studied a class of vertex
algebra-like structures, called $H_{D}$-quantum vertex algebras. The
notion of $H_{D}$-quantum vertex algebra generalizes
Etingof-Kazhdan's notion of braided vertex operator algebra
\cite{ek} in certain directions. In particular, the underlying space
of an $H_{D}$-quantum vertex algebra is a topologically free
$\C[[t]]$-module and the vertex operator map $Y$ is
$\C[[t]]$-linear, where the variable $t$ plays the same role as
$\hbar$ does in \cite{ek}. One of the generalizations is that the
braiding operator $\S$ (for the vertex operators on the algebra) is
allowed to have two (independent) spectral parameters, instead of
one. A fact is that general $H_{D}$-quantum vertex algebras fail to
satisfy the associativity for ordinary vertex algebras, though they
do satisfy a braided associativity, just as Etingof-Kazhdan's
braided vertex operator algebras do. (On the other hand,
Etingof-Kazhdan's quantum vertex operator algebras by definition
satisfy the associativity.)

We would like to thank Professor Hiromichi Yamada for organizing
this conference and for financial support.

\section{Nonlocal vertex algebras and their modules}
In this section, we review the basics on nonlocal vertex algebras
and their modules and quasi modules, and we give a summary of the
conceptual construction of nonlocal vertex algebras and their
(quasi) modules.

For this paper, letters such as $t, x,y, z, x_{0},x_{1},
x_{2},\dots$ are mutually commuting independent formal variables. We
shall use the formal variable notations and conventions as
established in \cite{flm} and \cite{fhl} (cf. \cite{ll}).  For this
paper we shall be working on the field $\C$ of complex numbers. For
any positive integer $r$, denote by $\C[[x_{1},x_{2},\dots,x_{r}]]$
the algebra of formal nonnegative power series and by
$\C((x_{1},\dots,x_{r}))$ the algebra of formal Laurent series which
are globally truncated with respect to all the variables. Note that
in the case $r=1$, $\C((x))$ is in fact a field. By
$\C_{*}(x_{1},x_{2},\dots,x_{r})$ we denote the extension of
$\C[[x_{1},x_{2},\dots,x_{r}]]$ by inverting all the nonzero
polynomials.

For any permutation $(i_{1},i_{2},\dots,i_{r})$ on $\{ 1,\dots,r\}$,
$\C((x_{i_{1}}))\cdots ((x_{i_{r}}))$ is a field, containing
$\C[[x_{1},\dots,x_{r}]]$ as a subalgebra, so (by a basic fact in
classical ring theory), there exists a unique algebra embedding
\begin{eqnarray}
\iota_{x_{i_{1}},\dots,x_{i_{r}}}: \C_{*}(x_{1},x_{2},\dots,x_{r})
\rightarrow \C((x_{i_{1}}))\cdots ((x_{i_{r}})),
\end{eqnarray}
extending the identity endomorphism of $\C[[x_{1},\dots,x_{r}]]$
(cf. \cite{fhl}). Note that both $\C_{*}(x_{1},\dots,x_{r})$ and
$\C((x_{i_{1}}))\cdots ((x_{i_{r}}))$ contain
$\C((x_{1},\dots,x_{r}))$ as a subalgebra. We see that
$\iota_{x_{i_{1}},\dots,x_{i_{r}}}$ preserves
$\C((x_{1},\dots,x_{r}))$ element-wise and that
$\iota_{x_{i_{1}},\dots,x_{i_{r}}}$ is also
$\C((x_{1},\dots,x_{r}))$-linear. For any nonzero polynomial $p\in
\C[x_{1},\dots,x_{r}]$, (as $\iota_{x_{i_{1}},\dots,x_{i_{r}}}$ is
an algebra homomorphism) we have
\begin{eqnarray}
\iota_{x_{i_{1}},\dots,x_{i_{r}}}(1/p)=p^{-1},
\end{eqnarray}
where $p^{-1}$ denotes the inverse of $p$ in $\C((x_{i_{1}}))\cdots
((x_{i_{r}}))$.

In the general field of vertex algebras, a very basic notion is that
of nonlocal vertex algebra, which generalizes the notion of vertex
algebra in the way that the notion of associative algebra
generalizes that of commutative associative algebra.

\bd{dnonlocalva} {\em  A {\em nonlocal vertex algebra} over $\C$ is
a vector space $V$, equipped with a linear map
\begin{eqnarray*}
Y: &&V\rightarrow \Hom (V,V((x)))\subset (\End V)[[x,x^{-1}]],\\
&&v\mapsto Y(v,x)
\end{eqnarray*}
and a vector ${\bf 1}\in V$, satisfying the conditions that $Y({\bf
1},x)=1$,
$$Y(v,x){\bf 1}\in V[[x]]\ \ \mbox{and }\ \lim_{x\rightarrow 0}Y(v,x){\bf 1}=v
\ \ \mbox{ for }v\in V,$$ and that for $u,v,w\in V$, there exists a
nonnegative integer $l$ such that
\begin{eqnarray}
(x_{0}+x_{2})^{l}Y(u,x_{0}+x_{2})Y(v,x_{2})w
=(x_{0}+x_{2})^{l}Y(Y(u,x_{0})v,x_{2})w.
\end{eqnarray}} \ed

\br{rnlva-field} {\em The notion of nonlocal vertex algebra, which
was defined in \cite{li-qva1}, is exactly the notion of weak
axiomatic $G_{1}$-vertex algebra in \cite{li-g1}, and it is
essentially the same as the notion of field algebra, studied in
\cite{bk} (cf. \cite{kac}).} \er

\bd{dmodule} {\em Let $V$ be a nonlocal vertex algebra. A {\em
$V$-module} is a vector space $W$, equipped with a linear map
\begin{eqnarray*}
Y_{W}:&& V\rightarrow \Hom (W,W((x)))\subset (\End
W)[[x,x^{-1}]],\\
&&v\mapsto Y_{W}(v,x),
\end{eqnarray*}
satisfying the conditions that
$$Y_{W}({\bf 1},x)=1_{W}\ \ \mbox{
(the identity operator on $W$)}$$
 and that for $u,v\in V,\; w\in W$,
there exists a nonnegative integer $l$ such that
$$(x_{0}+x_{2})^{l}Y_{W}(u,x_{0}+x_{2})Y_{W}(v,x_{2})w
=(x_{0}+x_{2})^{l}Y_{W}(Y(u,x_{0})v,x_{2})w.$$ A {\em quasi
$V$-module} is defined by simply replacing the last condition with
that for $u,v\in V,\; w\in W$, there exists a nonzero polynomial
$p(x_{1},x_{2})\in \C[x_{1},x_{2}]$ such that
\begin{eqnarray*}
p(x_{0}+x_{2},x_{2})Y_{W}(u,x_{0}+x_{2})Y_{W}(v,x_{2})w
=p(x_{0}+x_{2},x_{2})Y_{W}(Y(u,x_{0})v,x_{2})w.
\end{eqnarray*}}
\ed

\br{rinformation1} {\em The notion of module for a nonlocal vertex
algebra was introduced and a conceptual construction of nonlocal
vertex algebras and their modules was established in \cite{li-g1}.
The notion of quasi module for a vertex algebra was first introduced
and studied in \cite{li-gamma}, in order to associate vertex
algebras to a certain type of Lie algebras. Later, quasi modules for
nonlocal vertex algebras were studied in \cite{li-qva1} and a
general construction of nonlocal vertex algebras and their quasi
modules was established therein.} \er

Let $W$ be a general vector space over $\C$. Set
\begin{eqnarray}
\E(W)=\Hom (W,W((x)))\subset (\End W)[[x,x^{-1}]].
\end{eqnarray}
The identity operator on $W$, denoted by $1_{W}$, is a special
element of $\E(W)$.

\bd{dqcompatible} {\em A finite sequence $a_{1}(x),\dots,a_{r}(x)$
in $\E(W)$ is said to be {\em quasi compatible} if there exists a
nonzero polynomial $p(x,y)\in \C[x,y]$ such that
\begin{eqnarray}
\left(\prod_{1\le i<j\le r}p(x_{i},x_{j})\right) a_{1}(x_{1})\cdots
a_{r}(x_{r})\in \Hom (W,W((x_{1},\dots,x_{r}))).
\end{eqnarray}
The sequence $a_{1}(x),\dots,a_{r}(x)$
 is said to be {\em compatible} if there exists a nonnegative integer
$k$ such that
\begin{eqnarray}
\left(\prod_{1\le i<j\le r}(x_{i}-x_{j})^{k}\right)
a_{1}(x_{1})\cdots a_{r}(x_{r})\in \Hom (W,W((x_{1},\dots,x_{r}))).
\end{eqnarray}
Furthermore, a subset $T$ of $\E(W)$ is said to be {\em quasi
compatible} ({\em compatible})  if every finite sequence in $T$ is
quasi compatible (compatible).} \ed

Let $(a(x),b(x))$ be a quasi compatible ordered pair in $\E(W)$.
That is, there is a nonzero polynomial $p(x,y)\in \C[x,y]$ such that
\begin{eqnarray}\label{ehalf-relation}
p(x_{1},x_{2})a(x_{1})b(x_{2})\in \Hom (W,W((x_{1},x_{2}))).
\end{eqnarray}
We define $Y_{\E}(a(x),x_{0})b(x)\in \E(W)((x_{0}))$ by
\begin{eqnarray}
Y_{\E}(a(x),x_{0})b(x)=\iota_{x,x_{0}}\left(\frac{1}{p(x+x_{0},x)}\right)
\left(p(x_{1},x)a(x_{1})b(x)\right)|_{x_{1}=x+x_{0}}
\end{eqnarray}
and we then define $a(x)_{n}b(x)\in \E(W)$ for $n\in \Z$ by
\begin{eqnarray}
Y_{\E}(a(x),x_{0})b(x)=\sum_{n\in \Z}a(x)_{n}b(x) x_{0}^{-n-1}.
\end{eqnarray}
One can show that this is well defined; the expression on the
right-hand side is independent of the choice of $p(x,y)$. In this
way we have defined partial operations $(a(x),b(x))\mapsto
a(x)_{n}b(x)$ for $n\in \Z$ on $\E(W)$. We say that a quasi
compatible $\C$-subspace $U$ of $\E(W)$ is {\em $Y_{\E}$-closed} if
\begin{eqnarray}
a(x)_{n}b(x)\in U\ \ \ \mbox{ for }a(x),b(x)\in U,\; n\in \Z.
\end{eqnarray}

The main results of  \cite{li-qva1} (cf. \cite{li-g1}) can be
summarized as follows:

\bt{tquasi-main} Let $W$ be a general vector space over $\C$. a) For
any $Y_{\E}$-closed (quasi) compatible subspace $V$ of $\E(W)$, that
contains $1_{W}$, $(V,Y_{\E},1_{W})$ carries the structure of a
nonlocal vertex algebra with $W$ as a (quasi) module where
$Y_{W}(\alpha(x),x_{0})=\alpha(x_{0})$ for $\alpha(x)\in V$. b) For
every (quasi) compatible subset $U$ of $\E(W)$, there exists a
unique smallest $Y_{\E}$-closed (quasi) compatible subspace $\<U\>$
that contains $U$ and $1_{W}$, and  $\<U\>$ is a nonlocal vertex
algebra with $U$ as a generating subset and with $W$ as a (quasi)
module. \et

\section{Quantum vertex $\C((t))$-algebras and their modules}
In this section, we review the basics on  nonlocal vertex
$\C((t))$-algebras, (weak) quantum vertex $\C((t))$-algebras and
their modules, and we summarize the conceptual construction of weak
quantum vertex $\C((t))$-algebras and their (quasi) modules.

\bd{dt-algebra} {\em A {\em nonlocal vertex $\C((t))$-algebra} is a
nonlocal vertex algebra $V$ over $\C$, equipped with an
$\C((t))$-module structure, such that
\begin{eqnarray}
Y(f(t)u,x)(g(t)v)=f(t+x)g(t)Y(u,x)v
\end{eqnarray}
for $f(t),g(t)\in \C((t)),\; u, v\in V$. } \ed

For any field $K$, we denote by ${\bf Vec}(K)$ the category of
vector spaces over $K$. Note that for any vector space $W$ (over
$\C$), $\Hom (W,W((x)))$, alternatively denoted as $\E(W)$, is
naturally a $\C((x))$-module.

\bd{dmodule} {\em Let $V$ be a nonlocal vertex $\C((t))$-algebra. A
{\em $V$-module in category ${\bf Vec}(\C)$} is a module $(W,Y_{W})$
for $V$ viewed as a nonlocal vertex algebra over $\C$, satisfying
the condition that
\begin{eqnarray}
Y_{W}(f(t)v,x)w=f(x)Y_{W}(v,x)w \ \ \ \mbox{ for }f(t)\in \C((t)),\;
v\in V,\; w\in W.
\end{eqnarray}
A notion of {\em quasi $V$-module in category ${\bf Vec}(\C)$} is
defined in the obvious way---with the word ``module'' replaced by
``quasi module'' in the two places.} \ed

The following is a conceptual result which was obtained in
\cite{li-tqva}:

\bt{tconcrete-general} Let $W$ be a vector space over $\C$ and let
$U$ be any (quasi) compatible subset of $\E(W)$. Then $\C((x))\<U\>$
is the smallest $Y_{\E}$-closed (quasi) compatible
$\C((x))$-submodule of $\E(W)$, that contains $U$ and $1_{W}$, where
$\<U\>$ denotes the smallest $Y_{\E}$-closed $\C$-subspace of
$\E(W)$, that contains $U$ and $1_{W}$. Furthermore,
$(\C((x))\<U\>,Y_{\E},1_{W})$ carries the structure of a nonlocal
vertex $\C((t))$-algebra, where $$f(t)a(x)=f(x)a(x)\ \ \ \mbox{ for
}f(t)\in \C((t)),\; a(x)\in \C((x))\<U\>,$$ and $(W,Y_{W})$ carries
the structure of a (quasi) $\C((x))\<U\>$-module in category ${\bf
Vec}(\C)$, where $Y_{W}(a(x),x_{0})=a(x_{0})$ for $a(x)\in
\C((x))\<U\>$. \et

The following notion of weak quantum vertex $\C((t))$-algebra
singles out an important family of nonlocal vertex
$\C((t))$-algebras:

\bd{dtqva} {\em A {\em weak quantum vertex $\C((t))$-algebra} is a
nonlocal vertex $\C((t))$-algebra $V$ satisfying the condition that
for any $u,v\in V$, there exist
$$u^{(i)},v^{(i)}\in V,\; f_{i}(x_{1},x_{2})\in \C_{*}(x_{1},x_{2})
\ \ \ (i=1,\dots,r)$$
 such that
\begin{eqnarray}\label{et-jacobi-wqva}
& &x_{0}^{-1}\delta\left(\frac{x_{1}-x_{2}}{x_{0}}\right)
Y(u,x_{1})Y(v,x_{2})\nonumber\\
& &\ \ \ \ - x_{0}^{-1}\delta\left(\frac{x_{2}-x_{1}}{-x_{0}}\right)
\sum_{i=1}^{r}\iota_{t,x_{2},x_{1}}(f_{i}(t+x_{1},t+x_{2}))
Y(v^{(i)},x_{2})Y(u^{(i)},x_{1})\nonumber\\
&=&x_{2}^{-1}\delta\left(\frac{x_{1}-x_{0}}{x_{2}}\right)
Y(Y(u,x_{0})v,x_{2}).
\end{eqnarray}} \ed

A refinement of Theorem \ref{tconcrete-general} is that if a quasi
compatible subset $U$ of $\E(W)$ is of a certain type,
$\C((x))\<U\>$ is a weak quantum vertex $\C((t))$-algebra.

\bd{dpseudo-local} {\em Let $W$ be a vector space over $\C$. A
subset $U$ of $\E(W)$ is said to be {\em quasi
$\S(x_{1},x_{2})$-local}  if for any $a(x),b(x)\in U$, there exist
finitely many
$$u^{(i)}(x),\ v^{(i)}(x)\in U,\ \; f_{i}(x_{1},x_{2})\in \C_{*}(x_{1},x_{2})
\ \ (i=1,\dots,r)$$ such that
\begin{eqnarray}\label{epseudo-local}
p(x_{1},x_{2})a(x_{1})b(x_{2})
=\sum_{i=1}^{r}p(x_{1},x_{2})\iota_{x_{2},x_{1}}(f_{i}(x_{1},x_{2}))
u^{(i)}(x_{2})v^{(i)}(x_{1})
\end{eqnarray}
for some nonzero polynomial $p(x_{1},x_{2})\in \C[x_{1},x_{2}]$,
depending on $a(x)$ and $b(x)$. We say that $U$ is {\em
$\S(x_{1},x_{2})$-local} if for any $a(x),b(x)\in U$, there exist
$$u^{(i)}(x),v^{(i)}(x)\in U,\; f_{i}(x_{1},x_{2})\in
\C_{*}(x_{1},x_{2})\;\; (i=1,\dots, r)$$ such that
\begin{eqnarray}
(x_{1}-x_{2})^{k}a(x_{1})b(x_{2})
=(x_{1}-x_{2})^{k}\sum_{i=1}^{r}\iota_{x_{2},x_{1}}(f_{i}(x_{1},x_{2}))
u^{(i)}(x_{2})v^{(i)}(x_{1})
\end{eqnarray}
for some nonnegative integer $k$. } \ed

It was proved in \cite{li-tqva} that quasi $\S(x_{1},x_{2})$-local
subsets of $\E(W)$ are quasi compatible while
$\S(x_{1},x_{2})$-local subsets  are compatible. We have the
following fundamental result (see \cite{li-tqva}):

\bt{tmain} Let $W$ be a vector space over $\C$ and let $U$ be any
(quasi) $\S(x_{1},x_{2})$-local subset of $\E(W)$. Denote by $\<U\>$
the nonlocal vertex algebra over $\C$ generated by $U$. Then
$\C((x))\<U\>$ is a weak quantum vertex $\C((t))$-algebra with
$$f(t)a(x)=f(x)a(x)\ \ \mbox{ for }f(t)\in \C((t)),\; a(x)\in \C((x))\<U\>,$$
and $(W,Y_{W})$ is a (quasi) module for $\C((x))\<U\>$ in category
${\bf Vec}(\C)$, where
$$Y_{W}(a(x),x_{0})=a(x_{0})\ \ \ \mbox{ for }
a(x)\in \C((x))\<U\>.$$ \et

For nonlocal vertex $\C((t))$-algebras, there is another category of
modules which are like the adjoint modules.

\bd{dtmodule} {\em Let $V$ be a nonlocal vertex $\C((t))$-algebra. A
{\em (quasi) $V$-module in category ${\bf Vec}(\C((t)))$} is a
$\C((t))$-module $W$ which is also a (quasi) module for $V$ viewed
as a nonlocal vertex algebra over $\C$ such that
\begin{eqnarray}
Y_{W}(f(t)v,x)(g(t)w)=f(t+x)g(t)Y_{W}(v,x)w
\end{eqnarray}
for $f(t),g(t)\in \C((t)),\; v\in V,\; w\in W$.} \ed

\bd{dt1module} {\em Let $W$ be a $\C((t))$-module and let $t_{1}$ be
a formal variable. We define a $\C((t_{1}))$-module structure on
$\E(W)$ by
\begin{eqnarray}
f(t_{1})a(x)=f(t+x)a(x)\ \ \ \mbox{ for }f(t_{1})\in \C((t_{1})),\;
a(x)\in \E(W).
\end{eqnarray}}
\ed

With these notions we have:

\bp{pcompatible-t} Let $W$ be a $\C((t))$-module and let $U$ be a
compatible subset of $\E(W)$. Denote by $\<U\>$ the nonlocal vertex
algebra over $\C$ generated by $U$. Then $\C((t_{1}))\<U\>$ is a
nonlocal vertex $\C((t_{1}))$-algebra, and $W$, viewed as a
$\C((t_{1}))$-module with $f(t_{1})\in \C((t_{1}))$ acting as
$f(t)$, is a module in category ${\bf Vec}(\C((t_{1})))$. \ep

The notion of quantum vertex $\C((t))$-algebra involves quantum
Yang-Baxter operators. Let $H$ be a vector space over $\C$. A {\em
quantum Yang-Baxter operator} with two spectral parameters on $H$ is
a linear map
$$\S(x_{1},x_{2}): H\otimes H\rightarrow H\otimes H\otimes
\C_{*}(x_{1},x_{2})$$ satisfying the quantum Yang-Baxter equation
\begin{eqnarray}
\S_{12}(x_{1},x_{2})\S_{13}(x_{1},x_{3})\S_{23}(x_{2},x_{3})=
\S_{23}(x_{2},x_{3})\S_{13}(x_{1},x_{3})\S_{12}(x_{1},x_{2}),
\end{eqnarray}
where $\S_{ij}(x_{i},x_{j})$ are linear maps from $H^{\otimes
3}\rightarrow H^{\otimes 3}\otimes \C_{*}(x_{i},x_{j})$, defined by
$\S_{12}(x,z)=\S(x,z)\otimes 1,\ \ \S_{23}(x,z)=1\otimes \S(x,z)$,
and
$$\S_{13}(x,z)=P_{23}(\S(x,z)\otimes 1)P_{23}.$$
Furthermore, $\S(x_{1},x_{2})$ is said to be {\em unitary} if
\begin{eqnarray}
\S_{21}(x_{2},x_{1})\S(x_{1},x_{2})=1,
\end{eqnarray}
where $\S_{21}(x_{2},x_{1})=P\S(x_{2},x_{1})P$ with $P$ the flip
operator on $H\otimes H$.

\bd{dtqva-s} {\em A {\em quantum vertex $\C((t))$-algebra} is a weak
quantum vertex $\C((t))$-algebra $V$ equipped with a $\C$-linear
unitary quantum Yang-Baxter operator $\S(x_{1},x_{2})$ on $V$,
satisfying the conditions that
\begin{eqnarray}
\S(x_{1},x_{2})(f(t)u\otimes
g(t)v)=f(x_{1})g(x_{2})\S(x_{1},x_{2})(u\otimes v)
\end{eqnarray}
for $f(t),g(t)\in \C((t)),\; u,v\in V$, and that for any $u,v\in V$,
(\ref{et-jacobi-wqva}) holds with
$$\S(x_{1},x_{2})(u\otimes v)
=\sum_{i=1}^{r}u^{(i)}\otimes v^{(i)}\otimes f_{i}(x_{1},x_{2}), $$
and that
\begin{eqnarray}
& &[\D\otimes 1,\S(x_{1},x_{2})]=-\frac{\partial}{\partial x_{1}}\S(x_{1},x_{2}),\\
& &\S(x_{1},x_{2})(Y(x)\otimes 1) =(Y(x)\otimes
1)\S_{23}(x_{1},x_{2})\S_{13}(x_{1}+x-t,x_{2}),
\end{eqnarray}
where $\D$ is the $\C$-linear operator on $V$, defined by
$\D(v)=v_{-2}{\bf 1}$ for $v\in V$.} \ed

In the study of quantum vertex operator algebras, Etingof-Kazhdan
\cite{ek} introduced a notion of non-degeneracy, which has played a
very important role. This notion is also a very important tool in
the study of quantum vertex algebras in \cite{li-qva1}. The
following is a version of non-degeneracy for nonlocal vertex
$\C((t))$-algebras (cf. \cite{ek}):

\bd{dnondegenerate-ta} {\em Let $V$ be a nonlocal vertex
$\C((t))$-algebra. Denote by $V^{\otimes n}$ the tensor product in
the category of $\C$-vector spaces and define $ V^{\otimes
n}\boxtimes \C_{*}(x_{1},\dots,x_{n})$ to be the quotient space of
$V^{\otimes n}\otimes \C_{*}(x_{1},\dots,x_{n}) $ by the relations
$$ f_{1}(t)v^{(1)}\otimes \cdots \otimes f_{n}(t)v^{(n)}\otimes f
=  v^{(1)}\otimes \cdots \otimes v^{(n)}\otimes f_{1}(x_{1})\cdots
f_{n}(x_{n})f
$$
for $f\in \C_{*}(x_{1},\dots,x_{n}),\; f_{i}(t)\in \C((t)),\;
v^{(i)}\in V$ $(i=1,\dots,n)$. For each positive integer $n$, define
a $\C$-linear map
$$Z_{n}: V^{\otimes n}\boxtimes  \C_{*}(x_{1},\dots,x_{n})
\rightarrow V((x_{1}))\cdots ((x_{n}))$$ by
\begin{eqnarray*}
Z_{n} (v^{(1)}\otimes \cdots \otimes v^{(n)}\otimes f)=
\iota_{t,x_{1},\dots,x_{n}}f(t+x_{1},\dots, t+x_{n})
Y(v^{(1)},x_{1})\cdots Y(v^{(n)},x_{n}){\bf 1}.
\end{eqnarray*}
We say that $V$ is {\em non-degenerate} if for every positive
integer $n$, $Z_{n}$ is injective.} \ed

With this notion we have (see \cite{li-tqva}, cf. \cite{ek},
Proposition 1.11):

\bt{tek} Let $V$ be a weak quantum vertex $\C((t))$-algebra. Assume
that $V$ is non-degenerate. Then there exists a $\C$-linear map
$$\S(x_{1},x_{2}): V\otimes V\rightarrow V\otimes V\otimes
\C_{*}(x_{1},x_{2}),$$ which is uniquely determined by the condition
that for $u,v,w\in V$,
\begin{eqnarray}
&&(x_{1}-x_{2})^{k}Y(v,x_{2})Y(u,x_{1})w\nonumber\\
&=& (x_{1}-x_{2})^{k}Y(x_{1})(1\otimes
Y(x_{2}))(\S(x_{1}+t,x_{2}+t)(u\otimes v)\otimes w)
\end{eqnarray}
for some nonnegative integer $k$, depending only on $u$ and $v$.
Furthermore, $\S(x_{1},x_{2})$ is a unitary quantum Yang-Baxter
operator on $V$, and $V$ equipped with $\S(x_{1},x_{2})$ is a
quantum vertex $\C((t))$-algebra. \et

The following is a general result on non-degeneracy (see
\cite{li-tqva}, cf. \cite{li-qva2}):

\bp{pnon-degenerate} Let $V$ be a nonlocal vertex $\C((t))$-algebra
such that $V$ as a $V$-module is irreducible with
$\End_{V}(V^{\mod})=\C((t))$. Then $V$ is non-degenerate. \ep

\section{Quantum affine algebras and weak quantum vertex
$\C((t))$-algebras}

In this section we give a summary of the association of quantum
affine algebras with weak quantum vertex $\C((t))$-algebras and
their quasi modules.

First, we follow \cite{fj} (cf. \cite{drqaa}) to present the quantum
affine algebras. Let $\g$ be a finite-dimensional simple Lie algebra
of rank $l$ of type $A$, $D$, or $E$ and let $A=(a_{ij})$ be the
Cartan matrix. Let $q$ be a nonzero complex number. For $1\le i,j\le
l$, set
\begin{eqnarray}
f_{ij}(x)=(q^{a_{ij}}x-1)/(x-q^{a_{ij}})\in \C(x).
\end{eqnarray}
Then we set
\begin{eqnarray}
g_{ij}(x)^{\pm 1}=\iota_{x,0} f_{ij}(x)^{\pm 1}\in \C[[x]],
\end{eqnarray}
where $\iota_{x,0} f_{ij}(x)^{\pm 1}$ are the formal Taylor series
expansions of $f_{ij}(x)^{\pm 1}$ at $0$. Let $\Z_{+}$ denote the
set of positive integers. The quantum affine algebra
$U_{q}(\hat{\g})$ is (isomorphic to) the associative algebra with
identity $1$ with generators
\begin{eqnarray}
X_{ik}^{\pm},\ \ \phi_{im},\ \  \psi_{in},\ \ \gamma^{1/2},\ \
\gamma^{-1/2}
\end{eqnarray}
for $1\le i\le l, \;k\in \Z,\; m\in -\Z_{+},\;n\in \Z_{+}$, where
$\gamma^{\pm 1/2}$ are central, satisfying the relations below,
written in terms of the following generating functions in a formal
variable $z$:
\begin{eqnarray}
X_{i}^{\pm}(z)=\sum_{k\in \Z}X_{ik}^{\pm}z^{-k},\ \ \ \
\phi_{i}(z)=\sum_{m\in -\Z_{+}}\phi_{im}z^{-m},\ \ \ \
\psi_{i}(z)=\sum_{n\in \Z_{+}}\psi_{in}z^{-n}.
\end{eqnarray}
The relations are
\begin{eqnarray*}
& &\gamma^{1/2}\gamma^{-1/2}=\gamma^{-1/2}\gamma^{1/2}=1,\\
& &\phi_{i0}\psi_{i0}=\psi_{i0}\phi_{i0}=1,\\
& &[\phi_{i}(z),\phi_{j}(w)]=0,\ \ \ \
[\psi_{i}(z),\psi_{j}(w)]=0,\\
& &\phi_{i}(z)\psi_{j}(w)\phi_{i}(z)^{-1}\psi_{j}(w)^{-1}
=g_{ij}(z/w\gamma)/g_{ij}(z\gamma/w),\\
& &\phi_{i}(z)X^{\pm}_{j}(w)\phi_{i}(z)^{-1}
=g_{ij}(z/w\gamma^{\pm 1/2})^{\pm 1}X^{\pm}_{j}(w),\\
& &\psi_{i}(z)X^{\pm}_{j}(w)\psi_{i}(z)^{-1}
=g_{ij}(w/z\gamma^{\pm 1/2})^{\mp 1}X^{\pm}_{j}(w),\\
& &(z-q^{\pm 4a_{ij}}w)X^{\pm}_{i}(z)X^{\pm}_{j}(w)
=(q^{\pm 4a_{ij}}z-w)X^{\pm}_{j}(w)X^{\pm}_{i}(z),\\
& &[X^{+}_{i}(z),X^{-}_{j}(w)]=\frac{\delta_{ij}}{q-q^{-1}}
\left(\delta\left(\frac{z}{w\gamma}\right)\psi_{i}(w\gamma^{1/2})
-\delta\left(\frac{z\gamma}{w}\right)\phi_{i}(z\gamma^{1/2})\right),
\ \ \ \
\end{eqnarray*}
and there is one more set of relations of Serre type.

A $U_{q}(\hat{\g})$-module $W$ is said to be {\em restricted} if
 for any $w\in W$, $X_{ik}^{\pm}w=0$ and
$\psi_{ik}w=0$ for $1\le i\le l$ and for $k$ sufficiently large. We
say $W$ is of {\em level} $\ell\in \C$ if $\gamma^{\pm 1/2}$ act on
$W$ as scalars $q^{\pm \ell/4}$. (Rigorously speaking, one needs to
choose a branch of $\log q$.) We have (cf. \cite{li-qva1},
Proposition 4.9):

\bp{pqaffine} Let $q$ and $\ell$ be complex numbers with $q\ne 0$
and let $W$ be a restricted $U_{q}(\hat{\g})$-module of level
$\ell$. Set
$$U_{W}=\{ \phi_{i}(x), \psi_{i}(x), X_{i}^{\pm}(x)
\;|\; 1\le i\le l\}.$$ Then $U_{W}$ is a quasi
$\S(x_{1},x_{2})$-local subset of $\E(W)$ and $\C((x))\<U_{W}\>$ is
a weak quantum vertex $\C((t))$-algebra with $W$ as a quasi module
in category ${\bf Vec}(\C)$, where $\<U_{W}\>$ denotes the nonlocal
vertex algebra over $\C$ generated by $U_{W}$. \ep

With Proposition \ref{pqaffine} on hand, the remaining problem is to
determine the weak quantum vertex $\C((t))$-algebras $V_{W}$
explicitly and to show that they are quantum vertex
$\C((t))$-algebras, sufficiently by establishing the non-degeneracy.
We expect that just as with vertex algebras associated with affine
Lie algebras (cf. \cite{ll}), these weak quantum vertex
$\C((t))$-algebras are vacuum modules for certain associative
algebras derived {}from quantum affine algebras.

\end{document}